\documentclass{article}  
\usepackage{amssymb}              
\usepackage{amsthm}
\usepackage{amsmath}
\usepackage{eucal}
\usepackage{verbatim}
\usepackage[dvips]{graphicx}
\usepackage{multirow}
\usepackage{fancyhdr}
\usepackage{color}
\usepackage{enumerate}
\usepackage{calrsfs}
\usepackage{fullpage}
\usepackage{amsmath}
\usepackage{eufrak}

\usepackage{hyperref}

\newtheorem{theorem}{Theorem}
\newtheorem{lemma}{Lemma}

\newtheorem{corollary}{Corollary}

\makeatletter
\newcommand{\leqnomode}{\tagsleft@true}
\newcommand{\reqnomode}{\tagsleft@false}
\makeatother

\def\({\begin{eqnarray}}
\def\){\end{eqnarray}}
\def\[{\begin{eqnarray*}}
\def\]{\end{eqnarray*}}
\def\part#1#2{\frac{\partial #1}{\partial #2}}

\def\R{\mathbb{R}}
\def\N{\mathbb{N}}
\def\d{\mathrm{d}}
\def\tot#1#2{\frac{\d #1}{\d #2}}
\def\eps{\varepsilon}

\def\L{\mathcal{L}}

\def\L{\mathcal{L}}

\def\P{\mathrm{P}}

\def\L{\mathcal{L}}
\def\upsi{{\underline\psi}}

\begin{document}

\title{Asymptotic consensus with transmission and reaction delay: an~overview}   
\author{Jan Haskovec\footnote{Computer, Electrical and Mathematical Sciences \& Engineering, King Abdullah University of Science and Technology, 23955 Thuwal, KSA.
jan.haskovec@kaust.edu.sa}}

\date{}

\maketitle

\begin{abstract}
The aim of this paper is to provide a systematic overview of results on asymptotic consensus
for the Hegselmann-Krause-type model with delay and discuss the corresponding analytical tools.
We explain that two types (sources) of delay -- transmission and reaction --
are justifiable from the modeling point of view.
We consider both classical and normalized communication weights.
Studying a toy model with two agents only,
we develop an intuitive insight into what type of dynamics we can expect
from the systems.
In particular, we stress that with transmission-type delay,
asymptotic consensus can be reached with any length of the delay
(i.e., without smallness assumptions).
In contrast, the systems with reaction-type delay can only reach asymptotic consensus
if the delay is sufficiently small.

We formulate four theorems that establish asymptotic consensus
in the following scenarios:
(1) transmission-type delay with classical communication weights,
(2) transmission-type delay with normalized communication weights,
(3) reaction-type delay with symmetric communication weights,
(4) reaction-type delay with non-symmetric communication weights.
We explain how the methods of proof depend on the particular scenario:
direct estimates for (1), convexity arguments for (2), Lyapunov functional for (3)
and generalized Gronwall-Halanay inequality for (4).
\end{abstract}
\vspace{3mm}

\vspace{3mm}


\section{Introduction}\label{sec:Intro}
The Hegselmann-Krause model \cite{HK} of opinion dynamics
is a generic model of self-organization with applications in
biology, physics, socio-economic sciences and engineering
\cite{Camazine, Castellano, Jadbabaie, Krugman, Naldi, Vicsek, Xu}.
It describes a scenario with
a group of $N\in\N$ agents, each of them having his/her own opinion
or assessment of a certain quantity, represented by the real vector
$x_i\in\R^d$, $i\in\{1,2,\dots,N\}$, with $d\geq 1$.
Each of the agents communicates with all others and
revises his/her own opinion based on a weighted average of
the other agents' opinions. If this is done continuously in time,
we have $x_i=x_i(t)$ for $t\geq 0$ and
\( \label{eq:HK}
   \dot x_i(t) = \sum_{j=1}^N \psi_{ij}(t) (x_j(t) - x_i(t)), \qquad i\in [N],
\)
where here and in the sequel we denote $[N] := \{1, 2, \ldots,N\}$.
The communication weights $\psi_{ij}=\psi_{ij}(t)$ measure the intensity
of the influence between agents depending on the dissimilarity of their opinions.
This intensity, also called \emph{bounded confidence}, is usually
assumed to decay with increasing discord.
In the classical setting \cite{HK} the communication weights are given by
\(   \label{psi:nonr0}
   \psi_{ij}(t) :=  \frac{1}{N-1} \psi(|x_j(t) - x_i(t)|),
\)
with the nonnegative bounded \emph{influence function} $\psi:[0,\infty)\to [0,\infty)$,
also called \emph{communication rate}.
We adopt the assumption that $\psi(s) \leq 1$ for all $s\geq 0$, which implies
\[  
   \sum_{j\neq i} \psi_{ij} \leq 1 \qquad\mbox{for all } i\in [N].
\]
Let us note that this assumption is taken without loss of generality, since,
as long as $\psi$ is a globally bounded function, it can be achieved
by an appropriate rescaling of time.

As pointed out in \cite{MT} for the second-order version of \eqref{eq:HK}, aka the Cucker-Smale system \cite{CS1, CS2},
the scaling by $1/N$ in \eqref{psi:nonr0} has the drawback that the dynamics of an agent is modified by the
total number of agents even if its dynamics is only significantly influenced by a few nearby
agents. Therefore, \cite{MT} proposed to normalize the communication weights relative to the influence
of all other agents, without involving explicit dependence on their number,
\(\label{psi:r0}
   \psi_{ij}(t) := \frac{\psi(|x_j(t) - x_i(t)|)}{\sum_{\ell\neq i} \psi(|x_\ell(t) - x_i(t)|)}.
\)
The normalization removes the aforementioned drawback, but still allows
for a meaningful passage to the mean field limit $N\to\infty$, see \cite{MT, MT-SIAMRev}.
It also induces the following property of the normalized weights \eqref{psi:r0},
\(  \label{psi:conv}
   \sum_{j\neq i} \psi_{ij} = 1 \qquad\mbox{for all } i=1,\cdots, N,
\)
so that the term on the right-hand side of \eqref{eq:HK} is a convex combination
of the vectors $x_j(t) - x_i(t)$.

The main question is whether the dynamics \eqref{eq:HK} of continuously evolving opinions
will tend to an (asymptotic) emergence of one or more \emph{opinion clusters} formed
by agents with (almost) identical opinions.
\emph{Global consensus} is the state where all agents have the same opinion, i.e.,
$x_i=x_j$ for all $i,j \in [N]$.
Observe that the classical communication weights \eqref{psi:nonr0}
possess the symmetry property $\psi_{ij}=\psi_{ji}$, which implies
the conservation of the mean opinion $\frac{1}{N} \sum_{i=1}^N x_i$
along the dynamics induced by \eqref{eq:HK}.
Consequently, if the system reaches a global consensus,
its value is determined by the mean opinion of the initial datum.
On the other hand, the normalization in \eqref{psi:r0}
destroys the symmetry of the communication weights
and the global consensus value, if reached, cannot be
easily inferred from the initial datum and can be seen
as an emergent property of the system.

Various aspects of the consensus behavior of various modifications of \eqref{eq:HK}
have been studied in, e.g., \cite{Bha, Blondel, Canuto, Carro, Mohajer, Moreau, MT, JM, Wang, Wedin}.
For many applications in biological and socio-economical systems or control problems (for instance, swarm robotics \cite{Hamman, E3B, Valentini}),
it is natural to include a time delay in the model \eqref{eq:HK}.
From the modeling point of view, it is reasonable to consider two sources of delay:

\begin{itemize}
\item
\textbf{transmission-type delay,}
which reflects the fact that each agent may receive information from its surroundings
with a certain time lag due to finite speed of information transmission and/or perception.
Therefore, agent $i$ with opinion $x_i(t)$ receives at time $t>0$ the information about the opinion of agent $j$
in the form $x_j(t-\tau)$. In general, the length of the delay $\tau$ may depend both on the state of the system
and on the external conditions, i.e., it may be different for each pair $(i,j)$, it may vary in time
or follow a certain distribution.
For simplicity of the exposition, we consider a globally constant delay $\tau>0$ in this paper.
Then, system \eqref{eq:HK} is replaced by
\( \label{eq:HKprop}
   \dot x_i(t) = 
     \sum_{j\neq i} \psi_{ij}(t) (x_j(t-\tau) - x_i(t)), \qquad i\in [N].
\)
The communication rates $\psi_{ij} = \psi_{ij}(t)$ depend on the
history of the trajectories in a suitable way, as discussed below.

\item
\textbf{Reaction-type delay,}
which takes into account the time needed for the agents to react to the information that they already
received and processed. In this case, the agent $i$'s reaction executed at time $t$
is based upon the state of the system, including agent $i$'s state, at time $t-\tau$.
Again, the delay length $\tau$ may, in general, depend on both internal and external factors,
however, we shall restrict to the simple case of a fixed, constant delay $\tau>0$.
Then, system \eqref{eq:HK} is replaced by
\( \label{eq:HKreact}
   \dot x_i(t) = 
     \sum_{j\neq i} \psi_{ij}(t) (x_j(t-\tau(t)) - x_i(t-\tau)), \qquad i\in [N],
\)
where again the communication rates $\psi_{ij} = \psi_{ij}(t)$ depend on the history of the trajectories.
\end{itemize}
In both cases the system is equipped with the initial datum
\(\label{IC}
   x_i(t) = x^0_i(t),\qquad i=1,\cdots,N, \quad t \in [-\tau,0],
\)
with prescribed continuous trajectories $x^0_i\in\mathcal{C}([-\tau,0])$, $ i=1,\cdots,N$.

The choice of a particular expression for the communication rates $\psi_{ij} = \psi_{ij}(t)$
in \eqref{eq:HKprop} and \eqref{eq:HKreact} should be made depending on modeling considerations.
Again, there are two generic types of communication rates:

\begin{itemize}
\item
\textbf{Classical communication rates,} inspired by \eqref{psi:nonr0},
may take the form
\( \label{psi:prop:class}
   \psi_{ij}(t) = \frac{1}{N-1} \psi(|x_j(t - \tau) - x_i(t)|)
\)
for the system \eqref{eq:HKprop} with transmission-type delay, or
\(   \label{psi:react:class}
   \psi_{ij}(t) = \frac{1}{N-1} \psi(|x_j(t - \tau) - x_i(t-\tau)|)
\)
for the reaction-type delay system \eqref{eq:HKreact}.
However, we shall not set down any particular expression for $\psi_{ij}$
in the sequel. Instead, we characterize the classical communication weights by the property
\(  \label{psi:delay:subconv}
   \sum_{j\neq i} \psi_{ij}(t) \leq 1 \qquad\mbox{for all } t\geq 0 \mbox{ and } i\in [N].
\)
Note that \eqref{psi:delay:subconv} is obviously verified by both the expressions above
(recal the generic assumption $\psi\leq 1$, that can be made without loss of generality).

\item
\textbf{Normalized communication rates,} inspired by \eqref{psi:r0},
may take the form
\(  \label{psi:prop:norm}
   \psi_{ij}(t) := \frac{\psi(|x_j(t-\tau) - x_i(t)|)}{\sum_{\ell\neq i} \psi(|x_\ell(t-\tau) - x_i(t)|)}
\)
for the system \eqref{eq:HKprop} with transmission-type delay, or
\(   \label{psi:react:norm}
   \psi_{ij}(t) := \frac{\psi(|x_j(t-\tau) - x_i(t-\tau)|)}{\sum_{\ell\neq i} \psi(|x_\ell(t-\tau) - x_i(t-\tau)|)}
\)
for the reaction-type delay system \eqref{eq:HKreact}.
Again, we do not fix any particular form in the sequel.
Instead, we characterize the normalized communication weights by the property
\(  \label{psi:delay:conv}
   \sum_{j\neq i} \psi_{ij}(t) = 1 \qquad\mbox{for all } t\geq 0 \mbox{ and } i\in [N].
\)
Of course, \eqref{psi:delay:conv} is a special case of \eqref{psi:delay:subconv}.
However, it turns the terms on the right-hand side of \eqref{eq:HKprop}
and, resp., \eqref{eq:HKreact}, into convex combinations of the vectors $x_j(t-\tau) - x_i(t)$
and, resp., $x_j(t-\tau) - x_i(t-\tau)$. This peculiar convexity property facilitates
application of powerful analytical tools that provide stronger results
and/or more elegant proofs, as we shall demonstrate below.
\end{itemize}

Regardless of the particular form of the communication rates we shall consider,
we adopt the following \emph{global positivity assumption} for the whole paper:
If, for some $R>0$ and $T>0$, $|x_i(t)-x_j(t)| \leq R$ for all $ t\in [-\tau,T]$ and all $i, j \in [N]$,
then there exists $\upsi=\upsi(R)>0$ such that
\(  \label{ass:gp}
     \psi_{ij}(t) \geq \frac{\upsi}{N-1}  \quad \mbox{for all } t\in [-\tau,T] \mbox{ and all } i, j\in [N].
\)
Let us note that this assumption is satisfied by all the candidates
\eqref{psi:prop:class}--\eqref{psi:react:norm} proposed above
as soon as the influence function $\psi$ is globally positive.
Let us also note that the global positivity of $\psi$ is obviously
necessary for global consensus to be reached in general.
Indeed, if $\psi(r)=0$ for some $r>0$, then one would have steady states
formed by two clusters of particles at distance $r$ apart.

In order to guarantee existence of global $C^1$-solutions of
\eqref{eq:HKprop} and \eqref{eq:HKreact}, we adopt the generic assumption 
that all $\psi_{ij}=\psi_{ij}(t)$ are continuous functions of time.
This is verified by all the expressions \eqref{psi:prop:class}--\eqref{psi:react:norm}
as soon as $\psi\in\mathcal{C}([0,\infty))$.
Finally, let us stress that we do not impose any monotonicity property of $\psi$.

An important property that the consensus systems may or may not have is the \emph{symmetry}
and the resulting conservation of the mean
\(  \label{def:mean}
   \overline x(t) := \frac{1}{N} \sum_{i=1}^N x_i(t).
\)
Clearly, the transmission-type delay system \eqref{eq:HKprop} does not, in general,
conserve the mean even if one chooses symmetric communication rates
$\psi_{ij}=\psi_{ji}$. 
In contrast, the reaction-type system \eqref{eq:HKreact} conserves
$\overline x(t)$ as soon as the communication rates are symmetric,
like those given by \eqref{psi:react:class}. On the other hand, the normalized
weights \eqref{psi:react:norm} destroy the symmetry.
This is reflected in the methods that we shall use to prove asymptotic consensus
for the respective systems. For the symmetric case we shall apply an $\ell^2$-approach,
working with quadratic fluctuations of $x_i$ around the mean, while for the nonsymmetric case
an $\ell^\infty$-approach is more suitable. The conservation of the mean \eqref{def:mean}
is a significant property that allows us to obtain stronger results for the symmetric case.

The aim of this paper is to provide a systematic overview
of sufficient conditions for global asymptotic consensus
in the systems \eqref{eq:HKprop} and \eqref{eq:HKreact}
and the corresponding analytical tools.
Some of the results are somehow scattered in the recent literature and we feel that a unifying
perspective is beneficial. We also provide some new results and proofs.
In particular, we shall provide the following four theorems:
\begin{itemize}
\item
Theorem \ref{thm:Prop1} for the system \eqref{eq:HKprop} with communication rates satisfying \eqref{psi:delay:subconv},
\item
Theorem \ref{thm:Prop2} for \eqref{eq:HKprop} with normalized communication rates \eqref{psi:delay:conv},
\item
Theorem \ref{thm:React1} for the system  \eqref{eq:HKreact} with symmetric communication rates,
\item
Theorem \ref{thm:React2} for \eqref{eq:HKreact} with nonsymmetric communication rates \eqref{psi:delay:conv}.
\end{itemize}
The proofs are more or less closely following \cite{H:BLMS} and \cite{H:SIADS} for the system \eqref{eq:HKprop}.
The proof of Theorem \ref{thm:React1} is an adaptation of the method developed
in \cite{EHS} for a second-order system with noise.
The result for \eqref{eq:HKreact} with normalized (nonsymmetric) communication rates and its proof is
based on an adaptation of the method developed in \cite{H:JMAA}.

The paper is organized as follows. In Section \ref{sec:toy} we consider
the special case of \eqref{eq:HKprop} and \eqref{eq:HKreact} with two agents only, $N=2$,
and normalized communication weights.
In this setting the systems reduce to single linear delay equations
with known asymptotic behavior. They shall provide us intuition
on what type of qualitative behavior we may expect
from the solutions of the transmission-type delay consensus system \eqref{eq:HKprop}
in contrast to the reaction-type system \eqref{eq:HKreact}.
In Section \ref{sec:overview} we formulate the four Theorems about the asymptotic consensus
for the systems \eqref{eq:HKprop} and \eqref{eq:HKreact} with the two types
of communication rates discussed above.
We also explain the main ideas and methods of the proofs.
In Section \ref{sec:literature} we give an overview of previous literature on the topic
with brief comments on their results and methods.
In Section \ref{sec:aux} we prove an auxiliary Lemma of Gronwall-Halanay type
that shall be of later use.
Finally, in Sections \ref{sec:Prop} and \ref{sec:React} we provide proofs
of the four theorems formulated in Section \ref{sec:overview}.

\section{Toy models: $N=2$}\label{sec:toy}
To gain an intuition on what type of qualitative behavior we may expect
from the solutions of the delay consensus systems \eqref{eq:HKprop} and \eqref{eq:HKreact},
we consider the special case with two agents only, $N=2$, in the spatially one-dimensional setting $d=1$.
The goal of this section is to understand the different impact of
the transmission-type versus the reaction-type delay on the consensus dynamics.
Therefore, we consider agents with opinions $x_1(t)$, $x_2(t)$,
and, to simplify the situation even further, we choose the normalized communication weights.
Indeed, the convexity property \eqref{psi:delay:conv} immediately implies $\psi_{12}=\psi_{21}=1$.

With the above choices, the transmission-type delay system \eqref{eq:HKprop} reduces to
the linear system
\begin{align*} 
\begin{aligned}
   \dot x_1(t) &= x_2(t - \tau) - x_1(t), \\
   \dot x_2(t) &= x_1(t - \tau) - x_2(t).
\end{aligned}
\end{align*}
Defining 
$w:=x_1-x_2$, we have
\(
  \dot w(t) = -w(t - \tau) - w(t). \label{wEq}
\)
Assuming a solution of the form $w(t) = e^{\xi t}$ with complex $\xi\in\mathbb{C}$,
we obtain the characteristic equation
\[
   \xi  + e^{-\xi\tau} + 1 = 0.
\]
A simple inspection reveals that all roots $\xi$ have negative real part,
which implies that all solutions $w(t)$ of \eqref{wEq} tend to zero as $t\to\infty$.
Therefore, we have the asymptotic consensus $\lim_{t\to\infty} x_1(t)-x_2(t) = 0$,
for any value of the delay $\tau$.
This leads us to the intuitive expectation that {all} solutions of the general
system with transmission-type delay \eqref{eq:HKprop}
should converge to global consensus as $t\to\infty$, regardless of the length of the delay $\tau$.
Indeed, in Section \ref{sec:Prop} we shall prove this hypothesis for \eqref{eq:HKprop}
with both the classical communication rates characterized by \eqref{psi:delay:subconv}
and the normalized rates \eqref{psi:delay:conv}.

The situation is fundamentally different for the reaction-type delay in \eqref{eq:HKreact},
which for $N=2$ with \eqref{psi:delay:conv} turns into
\begin{align*} 
\begin{aligned}
   \dot x_1(t) &= x_2(t - \tau) - x_1(t-\tau), \\
   \dot x_2(t) &= x_1(t - \tau) - x_2(t-\tau).
\end{aligned}
\end{align*}
In this case we obtain the delay negative feedback equation for $w:=x_1-x_2$,
\(  \label{feedback}
    \dot w(t) = - 2w(t - \tau).
\)
An analysis of the corresponding characteristic equation
\[
   \xi + 2\tau e^{-\xi} = 0,
\]
with $\xi\in\mathbb{C}$, reveals that:
\begin{itemize}
\item
If $0 < 2\tau < e^{-1}$, then $u=0$ is asymptotically stable.
\item
If $e^{-1} < 2\tau < \pi/2$, then $u=0$ is asymptotically stable,
but every nontrivial solution of \eqref{feedback} is oscillatory.
\item
If $2\tau > \pi/2$, then $u=0$ is unstable.
\end{itemize}
In fact, if $2\tau < e^{-1}$, the solutions subject to the constant initial datum
never oscillate and tend to zero as $t\to\infty$.
If $2\tau$ becomes larger than $e^{-1}$ but
smaller than $\pi/2$, the nontrivial solutions must oscillate
(i.e., change sign infinitely many times as $t\to\infty$),
but the oscillations are damped and vanish as $t\to\infty$.
Finally, for $2\tau > \pi/2$ the nontrivial solutions oscillate with unbounded amplitude
as $t\to\infty$.
We refer to Chapter 2 of \cite{Smith} and \cite{Gyori-Ladas} for more details.

The different types of qualitative behavior of solutions of \eqref{wEq} versus \eqref{feedback}
can be intuitively understood by noting that the instantaneous negative feedback term $-w(t)$ in \eqref{wEq}
has a stabilizing effect, while the delay feedback term $-w(t-\tau)$ induces oscillations,
and this effect becomes stronger with longer delays.
Thus, the stabilizing effect of the instantaneous term $-w(t)$ in the right-hand side of \eqref{wEq}
is stronger than the destabilizing effect of the delay term $-w(t-\tau)$, regardless of the delay length $\tau$.
In \eqref{feedback} the stabilizing effect is not present, and therefore, for large enough delays, the solution
develops oscillations with divergent amplitude as $t\to\infty$. 
An analogous intuition applies to the general Hegselmann-Krause systems with transmission \eqref{eq:HKprop}
and reaction \eqref{eq:HKreact} type delay.
In \eqref{eq:HKprop} the instantaneous negative feedback term $-x_i(t)$ has a stabilizing effect,
which for any value of $\tau>0$ is stronger than the destabilization induced by $x_j(t-\tau(t))$.
On the other hand, in \eqref{eq:HKreact} the stabilization by $-x_i(t)$ does not take place,
so that one can expect asymptotic consensus only for short enough delays $\tau>0$.

Let us finally note that the traditional tool in study of asymptotic behavior of consensus-type models
is the quadratic fluctuation, which for $N=2$ reads simply $w(t)^2$.
Then, for \eqref{wEq} we have
\[
   \frac12 \tot{}{t} w(t)^2 = - w(t-\tau) w(t) - w(t)^2,
\]
and for \eqref{feedback},
\[
   \frac12 \tot{}{t} w(t)^2 = - 2 w(t-\tau) w(t).
\]
From here we clearly see the inherent analytical difficulty imposed by the presence of delay:
namely, that the term $w(t-\tau) w(t)$ does not have a definite sign as soon as oscillations occur.
Consequently, the quadratic fluctuation does not need to be a monotonically decaying function of time.
The goal of this paper is to present a systematic overview of methods that can be used to overcome
this difficulty for the Hegselmann-Krause models \eqref{eq:HKprop} and \eqref{eq:HKreact}.

\section{Discrete description: Overview of results and methods}\label{sec:overview}
 
 To formulate our main results, let us introduce the spatial diameter $d_x=d_x(t)$ of the agent group,
\(  \label{def:dx}
   d_x(t) := \max_{i,j \in [N]}|x_i(t) - x_j(t)|.
\)
The term \emph{global asymptotic consensus} for the systems \eqref{eq:HKprop} and \eqref{eq:HKreact} is then defined as the property
\( \label{def:consensus}
    \lim_{t\to\infty} d_x(t) = 0.
\)

\subsection{Transmission-type delay}
As announced in Section \ref{sec:toy}, we expect global asymptotic consensus to be achieved
for the system with transmission-type delay \eqref{eq:HKprop} for any $\tau>0$,
and for all initial data. Indeed, we have the following result:
 
 \begin{theorem}\label{thm:Prop1}
Let the communication weights $\psi_{ij}$ verify \eqref{psi:delay:subconv}
and the global positivity property \eqref{ass:gp}.
Then all solutions of \eqref{eq:HKprop} reach global asymptotic consensus in the sense of \eqref{def:consensus}.
\end{theorem}

We shall present a proof of Theorem \ref{thm:Prop1} in Section \ref{sec:Prop1},
based on the direct approach developed in \cite{H:BLMS}.
It is an $\ell^\infty$-type proof working with the group diameter $d_x=d_x(t)$
defined in \eqref{def:dx}.
We first consider the spatially one-dimensional setting
and derive a uniform bound on the group diameter in terms of the initial datum.
Then, we derive an explicit estimate on the diameter shrinkage
on finite time intervals. By an iterative argument we then conclude asymptotic convergence
of the group diameter to zero, i.e., consensus finding.
The disadvantage of this method is that it does not provide convergence rates.

An improved result can be obtained assuming that the communication weights satisfy the convexity property \eqref{psi:delay:conv},
e.g., the normalized communication weights \eqref{psi:prop:norm}.
In this case we have exponential convergence with calculable rate.
 
 \begin{theorem}\label{thm:Prop2}
Let the influence function be continuous and strictly positive on $[0, \infty)$
and let the weights $\psi_{ij}$ verify \eqref{psi:delay:conv}.
Then all solutions of \eqref{eq:HKprop} reach global asymptotic consensus as defined by \eqref{def:consensus}.
The decay of $d_x=d_x(t)$ to zero is exponential with calculable rate.
\end{theorem}

The proof of Theorem \ref{thm:Prop2}, developed in \cite{H:SIADS},
will be presented in Section \ref{sec:Prop2}.
It is again an $\ell^\infty$-type proof, and is significantly more elegant
than the rather technical proof of Theorem \ref{thm:Prop1} due to the benefit
of the convexity property \eqref{psi:delay:conv} of the communication weights.
It consists of three ingredients.
First, we prove nonexpansivity of the agent group, i.e., a uniform bound on the position radius of the agents.
Second, we establish a convexity argument, fundamentally based on \eqref{psi:delay:conv},
which provides a simple delay differential inequality.
Finally, we employ a Gronwall-Halanay-type inequality \cite{Halanay} to prove
that its solutions decay to zero exponentially, with rate
that can be calculated as a root of a nonlinear algebraic equation.

\subsection{Reaction-type delay}
Recalling the discussion in Section \ref{sec:toy}, we expect global asymptotic consensus to be achieved
only for small enough values of $\tau>0$ for the system \eqref{eq:HKreact} with reaction-type delay.
Consequently, we can only hope to derive a sufficient condition, in terms of $\tau>0$,
that guarantees asymptotic consensus for all solutions.
Moreover, due to the destabilizing effect of the reaction-type delay, as explained in Section \ref{sec:toy},
we need to adopt some additional/more restrictive assumptions
than in the case of transmission-type delay.

Our first result assumes that the communication weights
are symmetric, i.e., $\psi_{ij}(t) = \psi_{ji}(t)$ for all $i,j\in [N]$ and $t\geq 0$.
Moreover, we impose the upper bound
\(  \label{psi:N-1}
   \psi_{ij}(t) \leq \frac{1}{N-1},
\)
which is slightly more restrictive than \eqref{psi:delay:subconv}.
Note that both the symmetry condition and the upper bound \eqref{psi:N-1} are
verified by the classical weights \eqref{psi:react:class},
or, equally well by \eqref{psi:prop:class}, if $\psi\leq 1$,
which we assume without loss of generality.

\begin{theorem} \label{thm:React1}
Let the communication weights $\psi_{ij}$ verify \eqref{psi:N-1}, the global positivity \eqref{ass:gp}
and the symmetry property $\psi_{ij} = \psi_{ji}$.
If $\tau \leq 1/2$, then all solutions of \eqref{eq:HKreact} reach global asymptotic consensus as defined by \eqref{def:consensus}.
\end{theorem}

The proof of Theorem \ref{thm:React1} is an adaptation of the method developed in
\cite{EHS} and will be presented in Section \ref{sec:React1}.
It is an $\ell^2$-type proof, working with the quadratic fluctuation
of the agent positions, and is based on a construction of an appropriate Lyapunov functional.
It exploits the important consequence of the symmetry of the communication weigths,
namely, the conservation of the mean $\overline x(t) := \frac{1}{N} \sum_{i=1}^N x_i(t)$
along the solutions of \eqref{eq:HKreact}.
Indeed, we have
\(   \label{cons:mean}
   \tot{}{t} \overline x(t) = \frac{1}{N} \sum_{i=1}^N \sum_{j=1}^N \psi_{ij}(t) (x_j(t-\tau) - x_i(t-\tau)) = 0.
\)
The drawback of the proof is that, similarly as the proof of Theorem \ref{thm:Prop1},
it does not provide exponential convergence to consensus.

Our second result works for non-symmetric communication weights
satisfying \eqref{psi:delay:subconv}.
It provides exponential convergence to consensus with calculable rate,
however, at the expense of a tighter restriction on the length of $\tau>0$
and of the following assumption on the initial datum,
\(  \label{ICass}
     \max_{i\in [N]} |\dot x_i(t)| \leq d_x^0 \qquad\mbox{for almost all } t\in (-\tau,0),
\)
where we denoted
\(  \label{def:dx0}
   d_x^0 := \max_{t\in [-\tau,0]} d_x(t),
\)
with $d_x=d_x(t)$ defined in \eqref{def:dx}.
Let us note that \eqref{ICass} is trivially satisfied in the generic case
of constant initial datum on $[-\tau,0]$.
Moreover, we need to tighten the global positivity assumption \eqref{ass:gp}
in the sense that the lower bound on the communication weights $\psi_{ij}$
depends continuously on the diameter of the group $d_x$.
In particular, we assume that there exists a positive continuous function $\Psi\in\mathcal{C}([0,\infty))$ such that
\(  \label{ass:psi:cont}
   \psi_{ij}(t) \geq \frac{\Psi(d_x(t-\tau))}{N-1} \quad \mbox{for all } t\geq 0 \mbox{ and all } i, j\in [N].
\)
This is satisfied by the classical weights \eqref{psi:react:class}
as long as $\psi>0$, by choosing
\(  \label{uphi}
   \Psi(d_x) := \min_{s\in[0,d_x]} \psi(s).
\)
For the normalized weights \eqref{psi:react:norm} we have
\[
      \psi_{ij}(t) = \frac{\psi(|x_j(t-\tau) - x_i(t-\tau)|)}{\sum_{\ell\neq i} \psi(|x_\ell(t-\tau) - x_i(t-\tau)|)}
      \geq \frac{\psi(|x_j(t-\tau) - x_i(t-\tau)|)}{N-1},
\]
since, without loss of generality, $\psi\leq 1$. Consequently,
\eqref{ass:psi:cont} is again satisfied by choosing \eqref{uphi}.

\begin{theorem} \label{thm:React2}
Let the communication rates $\psi_{ij}$ satisfy \eqref{psi:delay:subconv}
and \eqref{ass:psi:cont}, and let the initial datum $x^0$ satisfy \eqref{ICass}.
Denoting $\upsi^0 := (N-1) \min_{i,j\in [N]}\psi_{ij} (0)$, assume that
\(     \label{ass:thm:4}
    4\tau < \upsi^0.
\)
Then there exists a unique $C\in(0,\upsi^0)$ such that
\[
   \upsi^0 - C = 4 e^{C\tau} \, \frac{e^{C\tau} - 1}{C},
\]
and the following estimate holds along solutions of \eqref{eq:HKreact}--\eqref{IC},
\(   \label{concl:thm:React2}
   d_x(t) \leq d_x^0 e^{-Ct} \qquad \mbox{for all } t \geq 0.
\)
\end{theorem}

The proof of Theorem \ref{thm:React2} is is an adaptation of the method developed
in \cite{H:JMAA} for a second-order system, aka the Cucker-Smale model.
It is based on a derivation of a delay integro-differential inequality
and an application of an appropriately generalized Gronwall-Halanay-type inequality.
We will present it in Section \ref{sec:React2}.

\section{Overview of previous literature}\label{sec:literature}

Exponential convergence to global consensus as $t\to\infty$ for the system \eqref{eq:HKprop}
with time varying delay $\tau=\tau(t)$
has been proved in \cite{CPP}, taking an $\ell^\infty$-approach and constructing a
Lyapunov functional that is a sum of the group diameter $d_x$ and an integral term.
The result holds under a set of conditions requiring smallness of the maximal time delay
in relation to the decay speed of the influence function and the fluctuation of the initial datum.
A variant of the proof for \eqref{eq:HKprop} with distributed time delay was given in \cite{Paolucci},
again a smallness assumption on the maximal delay.
These results are suboptimal since, as we state in Theorems \ref{thm:Prop1} and \ref{thm:Prop2},
the transmission-type delay system \eqref{eq:HKprop} reaches asymptotic consensus
without any restriction on $\tau$, and for all initial data.
Compared to Theorem \ref{thm:Prop1}, the advantage of \cite{CPP, Paolucci} is that it provides
exponential convergence. Theorem \ref{thm:Prop2} does
establish exponential convergence too, however, only for the convex interaction weights \eqref{psi:delay:conv}.

In \cite{Pignotti-Trelat} the authors studied convergence to consensus of the Cucker-Smale model
 with time-varying delays, which can be considered a second-order version of the system \eqref{eq:HKprop}.
 However, the authors imposed an assumption of a priori boundedness of the communication rates from below,
 which essentially reduces the model to first-order, i.e., a variant of \eqref{eq:HKprop}.
 The main result of \cite{Pignotti-Trelat} is convergence to consensus under a smallness assumption
 on the maximal time delay, which again is suboptimal in the context of transmission-type delay
 (cf. Theorem \ref{thm:Prop1}). The authors consider both symmetric interaction weights
 like \eqref{psi:nonr0}, where they apply an $\ell^2$-approach, and nonsymmetric weights
 like \eqref{psi:r0}, where they take an $\ell^\infty$-approach. In both cases the proof
 is based on a construction of a suitable Lyapunov functional.
 In \cite{Mauro} the optimal flocking result has been finally extablished for
 the Cucker-Smale model with transmission-type delay,
 where no smallness assumption on the delay has to be imposed.

The convexity argument (see Lemma \eqref{lem:geom} below)
that we shall apply for the case of normalized communication weights
satisfying \eqref{psi:delay:conv}
is an adaptation of the notion of coefficient of ergodicity introduced by Dobrushin \cite{Dobrushin}
in the context of Markov chains. 
It was later used to quantify the contractivity in first-order, discrete-time models of opinion dynamics (without delay)
in \cite{Krause00}. In \cite{HK-19} the class of \emph{scrambling matrices} was introduced, which can be characterized
as stochastic matrices of communication weights with sufficient connectivity between agents.
Then, the convexity argument, exploiting the stochasticity of the matrix,
is used to provide a bound on the contractivity of the system, which ultimately leads
to a proof of asymptotic consensus.
A similar argument was applied to second-order models (without delay) in \cite{MT-SIAMRev}.

Second-order variants of \eqref{eq:HKprop} and \eqref{eq:HKreact} were studied in \cite{Liu-Wu, EHS, HasMar, ChoiH1, ChoiH2}.
In particular, in \cite{Liu-Wu} and \cite{ChoiH1} the Cucker-Smale model with transmission-type delay
and normalized communication rates \eqref{psi:prop:norm} was considered
and asymptotic flocking was proved under a smallness condition on the delay,
related to the decay properties of the influence function $\psi$ and the velocity diameter of the initial datum.
In \cite{Liu-Wu} the authors study \eqref{eq:HKprop} with \eqref{psi:prop:norm}.
They generalize the technique of estimating ``maximal action of antisymmetric matrices",
developed in \cite{MT}, to the case with delay. However, this approach does not benefit from
the convexity property \eqref{psi:delay:conv} of the normalized communication weights and, therefore,
leads to suboptimal results. The method of proof in \cite{ChoiH1} is more direct,
based on deriving a system of dissipative delay differential inequalities and constructing
a Lyapunov functional.

Models of consensus finding with delay have been extensively studied
in the Engineering community, with results going back to \cite{E2}.
We provide only a short and very incomplete compilation of references here.
For linear problems, where the communication weights $\psi_{ij}$ are fixed,
stability criteria based on the frequency approach and on Lyapunov-Krasovskii techniques
were derived in \cite{E3A}. 
Proofs of delay-independent consensus and flocking
in nonlinear networks with multiple time-varying delays under mild assumptions
were provided in \cite{E4, E6}. These approaches are based on $\ell^\infty$-estimates obtained
using fundamental concepts from the non-negative matrix theory
and a Lyapunov-Krasovskii functional.
Consensus over directed static networks with arbitrary
finite communications delays was studied in \cite{Lu}.

\section{A Gronwall-Halanay-type lemma}\label{sec:aux}

We provide a proof of a Gronwall-Halanay-type lemma,
which establishes exponential decay of a nonnegative function, given that
it satisfies a certain integro-differential inequality.
The variant we present here is a slight
generalization of \cite[Lemma 2.5]{ChoiH1} and \cite[Lemma 3.3]{H:SIADS}.

 \begin{lemma}\label{lem:Halanay}
 Fix $\tau>0$ and let $u\in\mathcal{C}([-2\tau,\infty))$ be a nonnegative continuous function
with piecewise continuous derivative on $(0,\infty)$, such that
for almost all $t>0$ the differential inequality is satisfied,
\(   \label{Halanay:1}
   \tot{}{t} u(t) \leq \alpha \int_0^{\tau} u(t-s-\tau) \, \d\P(s) - \beta u(t),
\)
where $\P$ is a probability measure on $[0,\tau]$ and
$0 < \alpha < \beta$.
Then there exists $C\in(0,\beta-\alpha)$ such that
\(  \label{Halanay:2}
   \beta - C = \alpha \int_0^\tau e^{C (s + \tau)} \, \d\P(s),
\)
and the estimate holds
\( \label{Halanay:3}
   u(t) \leq \left( \max_{s\in[-2\tau,0]} u(s) \right) e^{-Ct} \qquad \mbox{for all } t \geq 0.
\)
\end{lemma}

\begin{proof}
It is easily seen that, for fixed $\tau>0$, the right-hand side of \eqref{Halanay:2}
is an unbounded increasing function of $C>0$, and that its limit as $C\to 0$ is $\alpha$.
On the other hand, the left-hand side of \eqref{Halanay:2} is decreasing
and equal to $\beta$ for $C=0$.
Consequently, \eqref{Halanay:2} is uniquely solvable whenever $\alpha < \beta$.
Then we have
\[
     \beta - C = \alpha \int_0^\tau e^{C (s + \tau)} \, \d\P(s) > \alpha,
\]
so that $C\in(0,\beta-\alpha)$.

Let us denote
\[   
    \bar u := \max_{s\in[-2\tau,0]} u(s), 
\]
and define the functions $v=v(t)$ and $w=w(t)$ as
\[
    v(t) :=  \left\{ \begin{array}{ll}
          \bar u & \textrm{for } t\leq 0,\\ [1mm]
           u(t) & \textrm{for } t >0,
  \end{array} \right.
  \qquad\qquad
    w(t) :=  \left\{ \begin{array}{ll}
          \bar u & \textrm{for } t\leq 0,\\ [1mm]
          \bar u e^{-Ct} & \textrm{for } t >0.
  \end{array} \right.
\]
Note that $v=v(t)$ verifies \eqref{Halanay:1} for almost all $t>0$.
For any fixed $\lambda > 1$ set
\[
   \mathcal{S}_\lambda := \Bigl\{ t \geq 0 : v(s) \leq \lambda w(s) \quad \mbox{for} \quad s \in [0,t] \Bigr\}.
\]
Since $0 \in \mathcal{S}_\lambda$, $T_\lambda := \sup \mathcal{S}_\lambda \geq 0$ exists. We claim that
\[
   T_\lambda = \infty \qquad \mbox{for any } \lambda >1.
\]
For contradiction, assume $T_\lambda < \infty$ for some $\lambda >1$.
Due to the continuity of $v=v(t)$ for $t>0$, there exists some $T_\lambda^* > T_\lambda$ such that
$v$ is differentiable at $T_\lambda^*$ and
\(  \label{est_derivatives}
   v(T_\lambda^*) > \lambda w(T_\lambda^*),\qquad   \tot{}{t} v(T_\lambda^*) > \lambda \tot{}{t} w(T_\lambda^*).
\)
Moreover, $T_\lambda^*$ can be chosen such that $T_\lambda^*-\tau < T_\lambda$.
Then we have, by construction, $v(s-\tau) \leq \lambda w(s-\tau)$ for all $s \leq T_\lambda^*$,
and from \eqref{Halanay:1} it follows
\[
   \tot{}{t} v(T_\lambda^*) &\leq& \alpha \int_0^\tau v(T_\lambda^*-s-\tau) \, \d\P(s) - \beta v(T_\lambda^*) \\
    &<& \alpha\lambda \int_0^\tau w(T_\lambda^*-s-\tau) \, \d\P(s) - \beta\lambda w(T_\lambda^*),
\]
where we also used \eqref{est_derivatives}.
Next, noting that, trivially, $w(t) \leq \bar u e^{-Ct}$ for all $t\in\R$, we have
\[
     \tot{}{t} v(T_\lambda^*) <  \alpha \lambda \bar u \int_0^\tau e^{-C (T_\lambda^* - s - \tau)} \, \d\P(s)  - \beta\lambda w(T_\lambda^*),
\]
and using the definition of $w(T_\lambda^*)$,
\[
    \tot{}{t} v(T_\lambda^*) <  
        \left( \alpha \int_0^\tau e^{C (s + \tau)} \, \d\P(s) - \beta \right) \lambda w(T_\lambda^*).
\]
Finally, assumption \eqref{Halanay:2} gives
\[
   \tot{}{t} v(T_\lambda^*)  < - C \lambda w(T_\lambda^*)  = \lambda \tot{}{t} w(T_\lambda^*),
\]
which is a contradiction to the second part of \eqref{est_derivatives}.

We conclude that, for every $\lambda>1$, $T_\lambda = \infty$ and $u(t) \leq \lambda w(t)$ for all $t\geq 0$.
Passing to the limit $\lambda \to 1$ yields the claim \eqref{Halanay:3}.
\end{proof}

We shall apply Lemma \ref{lem:Halanay} twice in this paper:
First in Section \ref{sec:Prop2}, where we choose 
the probability measure $\P$ to be a Dirac delta concentrated at zero.
We then use it again in Section \ref{sec:React2}, this time with $\P(s) := \tau^{-1} \chi_{[0,\tau]}(s)$,
where $\chi_{[0,\tau]}(s)$ denotes the characteristic function of the interval $[0,\tau]$.

\section{Asymptotic consensus with transmission-type delay}\label{sec:Prop}

A significant property of the system \eqref{eq:HKprop} with transmission-type delay
is that one can derive a uniform in time bound on the radius $R_x=R_x(t)$ of the agent group,
\(  \label{def:diam}
   R_x(t) := \max_{i\in [N]} |x_i(t)|,
\)
in terms of the radius of the initial datum, defined as
\[  
   R_x^0 := \max_{t\in [-\tau,0]} R_x(t).
\]
 In particular, we have the following result.

\begin{lemma}\label{lem:Rxbound}
Let the communication weights satisfy \eqref{psi:delay:subconv}.
Then, along the solutions of \eqref{eq:HKprop}, the diameter $R_x=R_x(t)$ defined in \eqref{def:diam} satisfies
\[
   R_x(t) \leq R_x^0 \qquad\mbox{for all } t\geq 0.
\]
\end{lemma}

This result, combined with the global positivity property \eqref{ass:gp}, directly implies
that the communication weights $\psi_{ij}=\psi_{ij}(t)$ are uniformly bounded below
by $\frac{\upsi}{N-1}$ with some $\upsi>0$.
We note that such a uniform a-priori bound is not available for the system \eqref{eq:HKreact} with reaction-type delay.
 
 The proof of Lemma \ref{lem:Rxbound} (or its variant for variable delay) was given in \cite{ChoiH1, CPP, H:SIADS}.
 Due to the fundamental role the result plays in the proof of asymptotic consensus,
 we provide it here for the sake of the reader.

\begin{proof}
Let us fix some $\eps >0$.
We shall prove that for all $t\geq 0$
\(  \label{R-eps}
   R_x(t) < R_x^0 + \eps.
\)
Obviously, $R_x(0) \leq R_x^0$, so that by continuity,
\eqref{R-eps} holds on the maximal interval $[0,T)$ for some $T>0$.
For contradiction, let us assume that $T<+\infty$.
Then we have
\(   \label{R-cont}
   \lim_{t\to T-}  R_x(t) = R_x^0 + \eps. 
\)
However, for any $i=1,\dots,N$, we have
\[
   \frac12 \tot{}{t} |x_i(t)|^2 &=&  \sum_{j\neq i} \psi_{ij}(t) \left[x_j(t-\tau)- x_i(t)\right]\cdot x_i(t) \\
     &=& \sum_{j\neq i}  \psi_{ij}(t) \left[ x_j(t-\tau)\cdot x_i(t) - |x_i(t)|^2\right].
\]
By definition, we have $|x_j(t-\tau)| < R_x^0 + \eps$ for all $j\in[N]$ and $t < T$,
so that with an application of the Cauchy-Schwarz inequality we arrive at
\[
   \frac12 \tot{}{t} |x_i(t)|^2 &\leq& \sum_{j\neq i} \psi_{ij}(t) \left[ \left(R_x^0 + \eps\right) |x_i(t)| - |x_i(t)|^2 \right]  \\
       &\leq& \left(R_x^0 + \eps\right) |x_i(t)| - |x_i(t)|^2,
\]
where we used \eqref{psi:delay:subconv} in the second line.
Now, if $|x_i(t)| \neq 0$, we use the identity $\frac12\tot{|x_i(t)|^2}{t} = |x_i(t)| \tot{|x_i(t)|}{t}$
and divide the above inequality by $|x_i(t)|$.
On the other hand, if $|x_i(t)| \equiv 0$ on an open subinterval of $(0,T)$,
then $\tot{|x_i(t)|}{t} \equiv 0 \leq R^0_x + \eps  - |x_i(t)|$ on this subinterval.
Thus, we obtain
\[
   \tot{}{t} |x_i(t)| \leq R^0_x  + \eps - |x_i(t)| \quad\mbox{for almost all } t\in(0,T),
\]
which implies
\[
   |x_i(t)| \leq \left(|x_i(0)| - (R^0_x + \eps) \right)e^{-t} + R^0_x + \eps  \qquad\mbox{for } t  < T.
\]
Consequently, with $|x_i(0)| \leq R_x^0$,
\[
   \lim_{t \to T-} \; \max_{1 \leq i \leq N} |x_i(t)| \leq -\eps e^{-T} + R_x^0 + \eps < R^0_x + \eps,
\]
which is a contradiction to \eqref{R-cont}.
We conclude that, indeed, $T = \infty$,
and complete the proof by taking the limit $\eps\to 0$.
\end{proof}

\subsection{Proof of Theorem \ref{thm:Prop1}}\label{sec:Prop1}
\def\I{\mathcal{I}}

The proof of Theorem \ref{thm:Prop1} was established in \cite{H:BLMS}
and we provide its main ideas here.
We first work in the spatially one-dimensional setting $d=1$, i.e., $x_i(t)\in\R$.
The 1D result is then easily generalized to the multi-dimensional
situation by considering arbitrary one-dimensional projections of \eqref{eq:HKprop}.
We remind the reader that in Theorem \ref{thm:Prop1} we only assumed the validity of
the uniform bound \eqref{psi:delay:subconv} for
the communication weights, together with the global positivity property \eqref{ass:gp}.

In the one-dimensional setting, one can establish the following analogue of Lemma \ref{lem:Rxbound}.

\begin{lemma}\label{lem:stay}
Let $d=1$ and define
\(   \label{mM}
   m := \min_{i\in [N]} \min_{t\in [-\tau,0]} x^0_i(t),\qquad
   M := \max_{i\in [N]} \max_{t\in [-\tau,0]} x^0_i(t).
\)
Let the communication weights satisfy \eqref{psi:delay:subconv}.
Then, along the solutions of \eqref{eq:HKprop},
\[   
   m \leq x_i(t) \leq M
\]
for all $i\in [N]$ and all $t\geq 0$.
\end{lemma}

\begin{proof}
The proof is similar to the proof of Lemma \ref{lem:Rxbound}
and can be found in \cite[Lemma 1]{H:BLMS}.
\end{proof}

Due to the global positivity assumption \eqref{ass:gp},
Lemma \ref{lem:stay} implies that
there exists $\upsi>0$ such that
\(  \label{upsi}
   \psi_{ij}(t) \geq \frac{\upsi}{N-1} \qquad \mbox{for all } t \geq 0 \mbox{ and } i, j\in [N].
\)

The following result establishes the fundament of the proof of asymptotic consensus.
It shows that on finite time intervals the spatial diameter of the agent group
shrinks by an explicit multiplicative factor.

\begin{lemma}\label{lem:shrink}
Let the assumptions of Lemma \ref{lem:stay} be verified and, moreover, assume $0 < m \leq M$.
Then, along the solutions of \eqref{eq:HKprop}, we have for all $i\in [N]$ and $t\in [5\tau,6\tau]$,
\[   
   m + \frac{\Gamma}{2} (M-m) \leq x_i(t) \leq M - \frac{\Gamma}{2} (M-m),
\]
where
\[  
   \Gamma := \left(1 - e^{-\frac{\upsi\tau}{N-1}} \right)^2 (1-e^{-\sigma}) e^{-6\tau} \frac{\upsi}{N-1},
\]
with $\upsi$ given by \eqref{upsi} and
\[   
   \sigma:=\min \left\{ \tau, \frac{M-m}{2M} \right\}.
\]
\end{lemma}

The proof is rather technical, based on direct estimation of the individual trajectories of the agents.
We refer to \cite{H:BLMS} for its details and proceed directly with the proof of Theorem \ref{thm:Prop1}.
Let us remind that we still consider the one-dimensional setting $d=1$.
We apply Lemma \ref{lem:shrink} iteratively on time intervals of total length $7\tau$.
Let us for $k\in\N$ denote $\I_k := [(6k-1)\tau,6k\tau]$ and
\[  
   m_k := \min_{i\in [N]} \min_{t\in \I_k} x_i(t),\qquad
   M_k := \max_{i\in [N]} \max_{t\in \I_k} x_i(t).
\]
Clearly, $m_0=m$ and $M_0=M$, with $m$ and $M$ given by \eqref{mM}.
Moreover, let us introduce the notation $D_k := M_k-m_k$ and
\[ 
   \Gamma_k := \left(1 - e^{-\frac{\upsi\tau}{N-1}} \right)^2 (1-e^{-\sigma_k}) e^{-6\tau} \frac{\upsi}{N-1},
\]
with $\sigma_k := \min \left\{ \tau, \frac{M_k-m_k}{2M_k} \right\}$. Note that $\Gamma_k\in (0,1)$ as long as $M_k>m_k$.

An application of Lemma \ref{lem:shrink} gives
\[
    m_0 + \frac{\Gamma_0}{2} (M_0-m_0) \leq x_i(t) \leq M_0 - \frac{\Gamma_0}{2} (M_0-m_0) \qquad\mbox{for } t\in\I_1.
\]
We thus have
\[
   D_1 = M_1-m_1 \leq (1 - \Gamma_0) (M_0-m_0) = (1-\Gamma_0) D_0.
\]
Iterating Lemma \ref{lem:shrink} gives
\[
   D_{k+1} \leq (1-\Gamma_k) D_k 
      \qquad\mbox{for all } k\in\N.
\]
Due to Lemma \ref{lem:stay} we have $M_k\leq M$ for all $k\in\N$, so that
denoting $\overline\sigma(D) := \min \left\{ \tau, \frac{D}{2M} \right\}$,
we have $\sigma_k \geq \overline\sigma(D_k)$ for all $k\in\N$.
Then, with
\[  
   \overline \Gamma(D) := \left(1 - e^{-\frac{\upsi\tau}{N-1}} \right)^2 (1-e^{-\overline\sigma(D)}) e^{-6\tau} \frac{\upsi}{N-1},
\]
we have $\Gamma_k \geq \overline \Gamma(D_k)$ for all $k\in\N$ and
\[   
   D_{k+1} \leq \left(1- \overline\Gamma(D_k)\right) D_k. 
\]
Clearly, the sequence $\{D_k\}_{k\in\N}$ is a nonnegative decreasing sequence, and denoting its limit $D$,
the limit passage $k\to\infty$ in the above inequality gives $\overline\Gamma(D)\leq 0$,
which immediately implies $D=0$.
To conclude the proof of Theorem \ref{thm:Prop1} in the one-dimensional setting,
we apply the uniform bound of Lemma \ref{lem:stay}, i.e., if the agent group is
contained in the spatial interval $[m_k, M_k]$ for $t\in\I_k$ 
then it remains contained in the same interval for all future times.

Finally, generalization of the proof to the spatially multi-dimensional setting
is facilitated by the observation that the claims of both Lemma \ref{lem:stay}
and Lemma \ref{lem:shrink} can be trivially adapted to projections
of the trajectories $x_i=x_i(t)$ to arbitrary one-dimensional subspaces of $\R^d$.
Indeed, for an arbitrary fixed vector $\xi\in\R^d$, we replace \eqref{eq:HKprop} with the projected system
\[
   \tot{}{t} (x_i(t)\cdot \xi) = \sum_{j\neq i} \psi_{ij}(t) (x_j(t-\tau) - x_i(t))\cdot\xi, \qquad i=1,\ldots,N,
\]
and by a simple adaptation of the above proofs we obtain
\[
   \lim_{t\to\infty} (x_i(t)-x_j(t))\cdot\xi = 0.
\]
for all $i,j\in[N].$ We conclude by choosing $\xi$'s as the basis vectors of $\R^d$.

\subsection{Proof of Theorem \ref{thm:Prop2}}\label{sec:Prop2}
The proof goes along the lines of \cite{H:SIADS} and is fundamentally
based on the convexity property \eqref{psi:delay:conv} of the communication weights
that is assumed in Theorem \ref{thm:Prop2}.
In particular, we have the following result, whose proof can be found in \cite[Lemma 3.2]{H:SIADS}.

\begin{lemma}\label{lem:geom}
Let $N\geq 3$ and $\{x_1, \dots, x_N\}\subset \R^d$ be any set of vectors in $\R^d$
and denote $d_x$ its diameter,
$$d_x:=\max_{1 \leq i,j \leq N}|x_i - x_j|.$$
Fix $i, k \in [N]$ such that $i\neq k$ and let
$\eta^i_j \geq 0$ for all $j\in[N] \setminus\{i\}$, 
and $\eta^k_j \geq 0$ for all $j\in[N] \setminus\{k\}$, 
such that
\[ 
   \sum_{j\neq i} \eta^i_j = 1,\qquad \sum_{j\neq k} \eta^k_j = 1.
\]
Let $\mu\geq 0$ be such that
\[ 
   \mu \geq \min \left\{ \min_{j\neq i} \eta^i_j,\; \min_{j\neq k} \eta^k_j \right\} \geq 0.
\]
Then
\[   
   \left| \sum_{j\neq i} \eta^i_j x_j - \sum_{j\neq k} \eta^k_j x_j \right| \leq (1-(N-2)\mu) d_x.
\]
\end{lemma}

In this section we shall apply Lemma \ref{lem:Halanay} with the probability measure $\P$
being a Dirac delta concentrated at zero. This gives the following result.

\begin{corollary}\label{corr:GH}
Let $u\in\mathcal{C}([-\tau,\infty))$ be a nonnegative continuous function
with piecewise continuous derivative on $(0,\infty)$, such that
for some constants $0 < \alpha < \beta$ the differential inequality is satisfied,
\[   
   \tot{}{t}  u(t) \leq \alpha u(t -\tau) - \beta u(t) \qquad\mbox{for almost all } t>0.
\]
Then there exists a unique solution $C\in (0,\beta-\alpha)$ of the equation
\begin{equation}\label{Halanay:1:2}
  \beta - C = \alpha e^{C\tau}
\end{equation}
and the estimate holds
\[
   u(t) \leq \left( \max_{s\in[-\tau,0]} u(s) \right) e^{-Ct} \qquad \mbox{for all } t \geq 0.
\]
\end{corollary}

\begin{proof}
Use Lemma \ref{lem:Halanay} with $\P(s) := \delta(s)$.
\end{proof}

We are now prepared to prove Theorem \ref{thm:Prop2}.

\begin{proof}

Due to the continuity of the trajectories $x_i=x_i(t)$,
there is an at most countable system of open, mutually disjoint
intervals $\{\mathcal{I}_\sigma\}_{\sigma\in\N}$ such that
$$
   \bigcup_{\sigma\in\N} \overline{\mathcal{I}_\sigma} = [0,\infty)
$$
and for each ${\sigma\in\N}$ there exist indices $i(\sigma)$, $k(\sigma)$
such that
$$
   d_x(t) = |x_{i(\sigma)}(t) - x_{k(\sigma)}(t)| \quad\mbox{for } t\in \mathcal{I}_\sigma.
$$
Then, using the abbreviated notation $i:=i(\sigma)$, $k:=k(\sigma)$,
we have for every $t\in \mathcal{I}_\sigma$,
\[
   \frac12 \tot{}{t} d_x(t)^2 &=& (\dot x_i(t) - \dot x_k(t))\cdot (x_i(t)-x_k(t))\\
      &=& 
       \left(\sum_{j\neq i} \psi_{ij}(t) (x_j(t-\tau) - x_i(t)) - \sum_{j\neq k} \psi_{kj}(t) (x_j(t-\tau) - x_k(t)) \right) \cdot (x_i(t)-x_k(t))\\
      &=&
      \left( \sum_{j\neq i} \psi_{ij}(t) x_j(t-\tau)  - \sum_{j\neq k} \psi_{kj}(t) x_j(t-\tau) \right) \cdot (x_i(t)-x_k(t)) - |x_i(t)-x_k(t)|^2,
\]
where we used the convexity property of the normalized weights \eqref{psi:conv}.
Lemma \ref{lem:Rxbound} together with the global positivity assumption \eqref{ass:gp}
provides some $\upsi\in (0,1)$ such that
\[
    \psi_{ij}(t) \geq \frac{\upsi}{N-1}  \qquad\mbox{for all } t\geq 0 \mbox{ and } i,j\in [N].
\]
We then apply Lemma \ref{lem:geom} with $\mu:=\frac{\underline{\psi}}{N-1}$, which gives
\[
   \left| \sum_{j\neq i} \psi_{ij} x_j(t-\tau)  - \sum_{j\neq k} \psi_{kj} x_j(t-\tau) \right| &\leq& (1-(N-2)\mu) d_x(t-\tau)  \\
     &=& \left(1- \frac{N-2}{N-1}\upsi\right) d_x(t-\tau) \,.
\]
Consequently, with the Cauchy-Schwartz inequality we have
\[
   \frac12 \tot{}{t} d_x(t)^2 \leq \left(1- \frac{N-2}{N-1}\upsi\right) d_x(t-\tau) d_x(t) - d_x(t)^2,
\]
which implies that for almost all $t>0$,
\[
   \tot{}{t} d_x(t) \leq \left(1- \frac{N-2}{N-1}\upsi\right) d_x(t-\tau) - d_x(t).
\]
An application of Corollary \ref{corr:GH} with $\alpha:= 1 - \frac{N-2}{N-1}\underline{\psi}$ and $\beta:=1$ gives then
the exponential decay
\[
   d_x(t) \leq \left( \max_{s\in[-\tau,0]} d_x(s) \right) e^{-Ct} \qquad \mbox{for } t \geq 0,
\]
where $C$ is the unique solution of \eqref{Halanay:1:2}, i.e.,
\[
    1 - C = \left(1- \frac{N-2}{N-1}\upsi\right) e^{C\tau}.
\]
\end{proof}

\section{Asymptotic consensus with reaction-type delay}\label{sec:React}
In this section we provide proofs of Theorems \ref{thm:React1} and \ref{thm:React2}
that give sufficient conditions for global asymptotic consensus 
in the system with reaction-type delay \eqref{eq:HKreact}.
Let us recall that, according to the discussion of the simplified model in Section \ref{sec:toy},
we cannot expect asymptotic consensus to take place for arbitrary delay lengths.
Indeed, the sufficient conditions of Theorems \ref{thm:React1} and \ref{thm:React2}
are formulated as smallness conditions on $\tau>0$.
This is related to the fact that, in contrast to the system with transmission-type delay (Lemma \ref{lem:Rxbound}),
no uniform bound on the radius of the agent group is available.
 
For the sake of legibility, in this section we shall use the shorthand notation $\widetilde x_j := x_j(t-\tau)$,
while $x_j$ means $x_j(t)$.

\subsection{Proof of Theorem \ref{thm:React1}} \label{sec:React1}
Let us recall that Theorem \ref{thm:React1} assumes symmetric interaction weights
\(   \label{sym}
   \psi_{ij}(t) = \psi_{ji}(t) \qquad \mbox{for all } i,j\in [N] \mbox{ and } t\geq 0,
\)
which implies the conservation of the mean \eqref{cons:mean} along the solutions of \eqref{eq:HKreact}.
Consequently, without loss of generality we may assume that
$$\overline x(t) = \frac{1}{N} \sum_{i=1}^N x_i(t) \equiv 0$$
for all $t\geq 0$.
We then introduce the quadratic fluctuation of the positions around the mean,
\(  \label{def:X}
   X(t):=\frac{1}{2(N-1)}\sum_{i=1}^N |x_{i}|^{2}.
\)
Moreover, for $t\geq 0$ we define the quantity
\(  \label{def:D}
   D(t):=\frac{1}{2(N-1)} \sum_{i=1}^N \sum_{j\neq i} {\psi}_{ij} |\widetilde{x}_{j} - \widetilde{x}_{i}|^{2}.
\)
We first derive an estimate on the dissipation of the quadratic fluctuation $X=X(t)$
in terms of $D=D(t)$. 

\begin{lemma} \label{lem:dX}
For any $\delta>0$ we have,
along the solutions of \eqref{eq:HKreact} with symmetric interaction weights \eqref{sym},
\begin{align} \label{dXest}
   \tot{}{t} X(t) \leq (\delta - 1) {D}(t) + {\tau}{\delta}^{-1} \int_{t-\tau}^{t} {D}(s) \,\d s
      \qquad\mbox{for all } t>\tau.
\end{align}
\end{lemma}

\begin{proof}
With \eqref{eq:HKreact} we have
\[
   \tot{}{t} X(t) &=& \frac{1}{N-1} \sum_{i=1}^N x_i \cdot \dot x_ i \\
     &=& \frac{1}{N-1} \sum_{i=1}^N \sum_{j\neq i} \psi_{ij}(t) (\widetilde x_j - \widetilde x_i) \cdot x_i \\
     &=& \frac{1}{N-1} \sum_{i=1}^N \sum_{j\neq i} \psi_{ij}(t) (\widetilde x_j - \widetilde x_i) \cdot \widetilde x_i
      + \frac{1}{N-1} \sum_{i=1}^N \sum_{j\neq i} \psi_{ij}(t) (\widetilde x_j - \widetilde x_i) \cdot (x_i - \widetilde x_i).
\]
For the first term of the right-hand side we apply the standard
symmetrization trick (exchange of summation indices $i\leftrightarrow j$, noting the symmetry of ${\psi}_{ij} = {\psi}_{ji}$),
\[
    \frac{1}{N-1} \sum_{i=1}^N \sum_{j=1}^N \psi_{ij}(t) (\widetilde x_j - \widetilde x_i) \cdot \widetilde x_i
      =  - \frac{1}{2(N-1)} \sum_{i=1}^N \sum_{j=1}^N \psi_{ij}(t) |\widetilde x_j - \widetilde x_i|^2 \,.
\]
Therefore, we have
\[
   \tot{}{t} X(t) = - \frac{1}{2(N-1)} \sum_{i=1}^N \sum_{j\neq i} \psi_{ij}(t) |\widetilde x_j - \widetilde x_i|^2
      + \frac{1}{N-1} \sum_{i=1}^N \sum_{j\neq i} \psi_{ij}(t) (\widetilde x_j - \widetilde x_i) \cdot (x_i - \widetilde x_i).
\]
For the last term we use the Young inequality with $\delta>0$,
\[
   \left|
      \frac{1}{N-1} \sum_{i=1}^N \sum_{j\neq i} \psi_{ij}(t) (\widetilde x_j - \widetilde x_i) \cdot (x_i - \widetilde x_i) \right|
     &\leq&
       \frac{\delta}{2(N-1)} \sum_{i=1}^N \sum_{j\neq i} \psi_{ij}(t) |\widetilde x_j - \widetilde x_i|^2
         + \frac{\delta^{-1}}{2(N-1)} \sum_{i=1}^N \sum_{j\neq i} \psi_{ij}(t) |x_i - \widetilde x_i|^2  \\
     &\leq&      
       \frac{\delta}{2(N-1)} \sum_{i=1}^N \sum_{j\neq i} \psi_{ij}(t) |\widetilde x_j - \widetilde x_i|^2
         + \frac{\delta^{-1}}{2(N-1)} \sum_{i=1}^N |x_i - \widetilde x_i|^2,
\]
where we used \eqref{psi:delay:subconv} for the second inequality.
Hence,
\( \nonumber
    \tot{}{t} X(t) &\leq& \frac{\delta-1}{2(N-1)} \sum_{i=1}^N \sum_{j\neq i} \psi_{ij}(t) |\widetilde x_j - \widetilde x_i|^2
        + \frac{\delta^{-1}}{2(N-1)} \sum_{i=1}^N |x_i - \widetilde x_i|^2 \\
        &=& (\delta-1) D(t) + \frac{\delta^{-1}}{2(N-1)} \sum_{i=1}^N |x_i - \widetilde x_i|^2.
      \label{lem:dX:1}
\)
Next, for $t>\tau$ we use \eqref{eq:HKreact} to evaluate the difference
$x_{i}-\widetilde{x}_{i}$,
\[
   x_{i}-\widetilde{x}_{i} = \int_{t-\tau}^{t} \dot x_{i}(s)\, \d s
    = \sum_{j\neq i} \int_{t-\tau}^{t}  {\psi}_{ij}(s)({x}_{j}(s-\tau)-{x}_{i}(s-\tau))\, \d s.
\]
Taking the square and summing over $i\in [N]$ we arrive at
\[
   \sum_{i=1}^N |x_i - \widetilde{x}_{i}|^{2}
   &=& \sum_{i=1}^N \left|  \sum_{j\neq i} \int_{t-\tau}^{t}
          {\psi}_{ij}(s)({x}_{j}(s-\tau) - {x}_{i}(s-\tau)) \,\d s \right|^2  \\
   &\leq& (N-1) \sum_{i=1}^N\sum_{j\neq i} \left|
          \int_{t-\tau}^{t} {\psi}_{ij}(s)({x}_{j}(s-\tau) - {x}_{i}(s-\tau)) \,\d s \right|^{2} \\
   &\leq& \tau \sum_{i=1}^N\sum_{j\neq i}  \int_{t-\tau}^{t}
          {\psi}_{ij}(s) |{x}_{j}(s-\tau)-{x}_{i}(s-\tau)|^{2} \,\d s  \\
   &=& 2 (N-1) \tau \int_{t-\tau}^{t} {D}(s)\,\d s.
\]
The first inequality is the Cauchy-Schwartz inequality for the summation term,
i.e. $\left|\sum_{j\neq i} a_{i}\right|^2 \leq (N-1)\sum_{j\neq i} |a_{i}|^2$,
and the second Cauchy-Schwartz inequality for the integral term,
together with the bound $\psi_{ij}\leq 1/(N-1)$ provided by assumption \eqref{psi:N-1}.
Inserting into \eqref{lem:dX:1} finally gives \eqref{dXest}.
\end{proof}

We now provide the proof of Theorem \ref{thm:React1}.

\begin{proof}
For some $\lambda>0$ we define the functional
\[   
   \L(t) := X(t) + \lambda \int_{t-\tau}^{t}\int_{\theta}^{t} {D}(s) \, \d s \, \d\theta,
\]
where $X=X(t)$ is the quadratic velocity fluctuation \eqref{def:X}
and $D=D(t)$ defined in \eqref{def:D}.

The time derivative of the second term in $\L(t)$ yields
\begin{align*} 
   \tot{}{t} \left( \lambda \int_{t-\tau}^{t}\int_{\theta}^{t}  {D}(s) \, \d s \, \d\theta \right)
       = \lambda\tau {D}(t) - \lambda \int_{t-\tau}^{t} {D}(s) \, \d s.
\end{align*}
Combining this with \eqref{dXest}, with the  choice $\lambda:=\tau\delta^{-1}$ so that we eliminate the integral term,
we obtain
 \[
   \tot{}{t} \L(t) \leq \left(\delta -1 + \tau^2\delta^{-1} \right) {D}(t) \qquad\mbox{for all } t\geq \tau.
\]
Optimizing in $\delta>0$ leads to $\delta:=\tau$ and
\(   \label{LyapDecay}
   \tot{}{t} \L(t) \leq \left(2\tau -1 \right) {D}(t) \qquad\mbox{for all } t\geq \tau.
\)
Therefore, we have
\[
   X(t) \leq \L(t) \leq \L(\tau)  \qquad\mbox{for all } t\geq \tau,
\]
and, recalling the conservation of the mean $\bar x(t)$,
the uniform boundedness of $X(t)$ implies uniform boundedness of $|x_i-x_j|$
for all $i,j\in [N]$.
Then, the global positivity assumption \eqref{ass:gp} gives some $\upsi>0$ such that
\[
   \psi_{ij}(t) \geq \frac{\upsi}{N-1}  \qquad\mbox{for all } t\geq \tau.
\]
Then, with $\bar x(t)\equiv 0$, we have
\[
      D(t) = \frac{1}{2(N-1)} \sum_{i=1}^N \sum_{j\neq i} {\psi}_{ij} |\widetilde{x}_{j} - \widetilde{x}_{i}|^{2} 
         \geq \frac{\upsi}{2(N-1)^2} \sum_{i=1}^N \sum_{j\neq i} |\widetilde{x}_{j} - \widetilde{x}_{i}|^{2} 
         = 2 \upsi X(t-\tau).
\]
Integrating \eqref{LyapDecay} on $[\tau,t]$ gives
\[
   X(t) \leq
   \L(t) &\leq& \L(\tau) + \left(2\tau -1 \right) \int_\tau^t {D}(s) \, \d s \\
          &\leq&  \L(\tau) + 2 \left(2\tau -1 \right)\upsi \int_0^{t-\tau} {X}(s) \, \d s.
\]
Taking limes superior as $t\to\infty$ in both sides, we have
\(   \label{limsup}
     2 \left(1-2\tau \right)\upsi \int_0^{\infty} {X}(s) \, \d s + \limsup_{t\to\infty} X(t) \leq \L(\tau),
\)
so that $X=X(t)$ is integrable at infinity if $\tau<1/2$.
Barbalat's lemma \cite{Barbalat} then states that $X(t)\to 0$ as $t\to\infty$
as soon as $\tot{}{t} X(t)$ is uniformly bounded, say, for $t>\tau$.
This follows easily by an application of Lemma \ref{lem:dX} with $\delta:=1$,
\[
    \tot{}{t} X(t) \leq \tau \int_{t-\tau}^{t} {D}(s) \,\d s 
       \leq \tau^2 \max_{s\in [t-\tau,t]} D(s),
\]
combined with \eqref{def:D} and \eqref{psi:N-1},
\[
      D(t) &=& \frac{1}{2(N-1)} \sum_{i=1}^N \sum_{j\neq i} {\psi}_{ij} |\widetilde{x}_{j} - \widetilde{x}_{i}|^{2}  \\
         &\leq&  \frac{1}{(N-1)^2} \sum_{i=1}^N \sum_{j\neq i} |\widetilde{x}_{j}|^2 + |\widetilde{x}_{i}|^{2}  \\
         &\leq&  4 X(t-\tau)
\]
and the observation that \eqref{limsup} implies uniform boundedness of $X=X(t)$.
\end{proof}

\subsection{Proof of Theorem \ref{thm:React2}}\label{sec:React2}
Here we consider the system with reaction-type delay \eqref{eq:HKreact}
with non-symmetric communication weights.
The proof of asymptotic consensus presented in this section is new.
We assume that the communication weights satisfy the global upper bound \eqref{psi:delay:subconv}
and depend continuously on the diameter of the group in the sense of \eqref{uphi}.
We also impose the technical assumption \eqref{ICass} on the initial datum.
For the sake of more elegant expressions, let us extend the definition
of $d_x=d_x(t)$ to $t<0$ as follows,
\(  \label{redef:dx}
   d_x(t) := \begin{cases}
      \quad \max_{i,j \in [N]} |x_i(t) - x_j(t)| &\qquad\mbox{for } t > 0,\\
      \quad d_x^0 &\qquad\mbox{for } t \leq 0,
      \end{cases}
\)
where $d_x^0$ is given by \eqref{def:dx0}.

\begin{lemma} \label{lem:4}
Let the communication weights satisfy \eqref{psi:delay:subconv}
and assume that there exists $\upsi >0$ such that, for some $T>0$,
\(  \label{ass:lem:4}
   \psi_{ij}(t) \geq \frac{\upsi}{N-1} \qquad\mbox{for all } t\in [0,T] \mbox{ and } i,j\in [N].
\)
Moreover, let the initial datum verify \eqref{ICass}.
Then, along the solutions of \eqref{eq:HKreact}, 
\[  
   \tot{}{t} d_x (t) \leq 4 \int_{t-\tau}^t d_x(s-\tau) \d s - \upsi d_x (t)
\]
for almost all $t\in (0,T)$.
\end{lemma}

\begin{proof}
Due to the continuity of the trajectories $x_i=x_i(t)$,
there is an at most countable system of open, mutually disjoint
intervals $\{\mathcal{I}_\sigma\}_{\sigma\in\N}$ such that
$$
   \bigcup_{\sigma\in\N} \overline{\mathcal{I}_\sigma} = [0,\infty)
$$
and for each ${\sigma\in\N}$ there exist indices $i(\sigma)$, $k(\sigma)$
such that
$$
   d_x(t) = |x_{i(\sigma)}(t) - x_{k(\sigma)}(t)| \quad\mbox{for } t\in \mathcal{I}_\sigma.
$$
Then, using the abbreviated notation $i:=i(\sigma)$, $k:=k(\sigma)$,
we have for every $t\in \mathcal{I}_\sigma$,
\(
   \frac12 \tot{}{t} d_x(t)^2 &=& (\dot x_i - \dot x_k)\cdot (x_i-x_k)    \nonumber\\
      &=& 
       \left(\sum_{j\neq i} \psi_{ij}(t) (\widetilde x_j - \widetilde x_i) - \sum_{j\neq k} \psi_{kj}(t) (\widetilde x_j - \widetilde x_k) \right) \cdot (x_i-x_k).
       \label{4:2}
\)
We process the first term of the right-hand side as follows,
\(   \label{4:1}
    \sum_{j\neq i} \psi_{ij}(t) (\widetilde x_j - \widetilde x_i) \cdot (x_i-x_k) &=&
      \sum_{j\neq i} \psi_{ij}(t) (\widetilde x_j - x_j + x_i - \widetilde x_i) \cdot (x_i-x_k) \\
         && + \sum_{j\neq i} \psi_{ij}(t) (x_j- x_i) \cdot (x_i-x_k).  \nonumber
\)
We estimate the difference $|\widetilde x_j - x_j|$ by
\[
   |\widetilde x_j - x_j| &\leq& \int_{t-\tau}^t |\dot x_j(s)| \,\d s \\
      &\leq&   \int_{-[t-\tau]^-}^0 |\dot x_j(s)| \,\d s  + \int_{[t-\tau]^+}^t \sum_{\ell\neq j} \psi_{j\ell}(s) |x_\ell(s-\tau) - x_j(s-\tau)| \,\d s \\
      &\leq&   \int_{-[t-\tau]^-}^0 d_x^0 \,\d s  + \int_{[t-\tau]^+}^t \sum_{\ell\neq j} \psi_{j\ell}(s) \, d_x(s-\tau) \,\d s \\
      &\leq& \int_{t-\tau}^t d_x(s-\tau) \, \d s,
\]
where for the third inequality we used the assumption \eqref{ICass} and for the 
last inequality we used the property \eqref{psi:delay:subconv},
together with \eqref{redef:dx}.
Performing an analogous estimate for the term $|\widetilde x_i - x_i|$
and using the Cauchy-Schwartz inequality and, again,  \eqref{psi:delay:subconv}, we arrive at
\[
   \sum_{j\neq i} \psi_{ij}(t) (\widetilde x_j - x_j + x_i - \widetilde x_i) \cdot (x_i-x_k) &\leq&
      \sum_{j\neq i} \psi_{ij}(t) \left( |\widetilde x_j - x_j| + |\widetilde x_i - x_i| \right) |x_i-x_k| \\
      &\leq&
      2\, d_x(t) \int_{t-\tau}^t d_x(s-\tau) \,\d s.
\]
To estimate the second term of the right-hand side of \eqref{4:1}, observe that, using the Cauchy-Schwarz inequality, we have
\[
   (x_j- x_i)\cdot (x_i-x_k) &=& (x_j-x_k)\cdot(x_i-x_k) - |x_i-x_k|^2 \\
      &\leq& |x_i-x_k| \bigl( |x_j-x_k| - |x_i-x_k| \bigr) \leq 0,
\]
since, by definition, $|x_j-x_k| \leq d_x = |x_i-x_k|$.
Then, by assumption \eqref{ass:lem:4}, we have
\[
   \sum_{j\neq i} \psi_{ij}(t) (x_j- x_i) \cdot (x_i-x_k) \leq
        \frac{\upsi}{N-1} \sum_{j=1}^N (x_j- x_i) \cdot (x_i-x_k).
\]
Repeating the same steps for the second term of the right-hand side of \eqref{4:2}, we finally arrive at
\[
      \frac12 \tot{}{t} d_x(t)^2 &\leq&
           4\, d_x(t) \int_{t-\tau}^t d_x(s-\tau) \,\d s
         + \frac{\upsi}{N-1} \left( \sum_{j=1}^N (x_j- x_i) \cdot (x_i-x_k) - \sum_{j=1}^N (x_j- x_k) \cdot (x_i-x_k) \right) \\
         &=&
           4\, d_x(t) \int_{t-\tau}^t d_x(s-\tau) \,\d s - \frac{N \upsi}{N-1} |x_i-x_k|^2 \\
       &\leq&
           4\, d_x(t) \int_{t-\tau}^t d_x(s-\tau) \,\d s - {\upsi}\, d_x(t)^2.
\]
Consequently, for almost all $t\in (0,T)$,
\[
      \tot{}{t} d_x(t) \leq  4 \int_{t-\tau}^t d_x(s-\tau) \,\d s - {\upsi}\, d_x(t).
\]
\end{proof}

In the sequel we shall apply the Gronwall-Halanay Lemma \ref{lem:Halanay} with the probability measure
$\P(s) := \tau^{-1} \chi_{[0,\tau]}(s)$,
where $\chi_{[0,\tau]}(s)$ is the characteristic function of the interval $[0,\tau]$.
This gives the following result.

\begin{corollary}\label{cor:Halanay:2}
Let $u\in\mathcal{C}([-2\tau,\infty))$ be a nonnegative continuous function
with piecewise continuous derivative on $(0,\infty)$, such that
for almost all $t>0$ the differential inequality is satisfied,
\[  
   \tot{}{t} u(t) \leq \frac{\alpha}{\tau} \int_{t-\tau}^t u(s-\tau) \, \d s - \beta u(t).
\]
Moreover, let $0 < \alpha < \beta$.
Then there exists $C\in(0,\beta-\alpha)$ such that
\(  \label{Halanay:2:2}
   \beta - C = \alpha e^{C\tau} \, \frac{e^{C\tau} - 1}{C\tau},
\)
and the estimate holds
\[  
   u(t) \leq \left( \max_{s\in[-2\tau,0]} u(s) \right) e^{-Ct} \qquad \mbox{for all } t \geq 0.
\]
\end{corollary}

\begin{proof}
Apply Lemma \ref{lem:Halanay} with $\P(s) := \tau^{-1} \chi_{[0,\tau]}(s)$, noting that
\[
    \int_0^\tau e^{C (s + \tau)} \, \d\P(s) = \frac{1}{\tau} \int_0^\tau e^{C (s + \tau)} \d s 
       = e^{C\tau} \, \frac{e^{C\tau} - 1}{C\tau}.
\]
\end{proof}

We are now ready to carry out the proof of Theorem \ref{thm:React2}.

\begin{proof}
We define the set
\[
   \mathcal{S} := \bigl\{ T\geq 0; \, d_x(t) \leq d_x^0 \mbox{ for all } t\in [0,T] \bigr\}.
\]
Obviously, $0\in S$, so that there exists $T^* := \sup \mathcal{S}$.
We claim that $T^* = \infty$.
For contradiction, assume that $T^* < \infty$, then we have
$d_x(t) \leq d_x^0$ for all $t\in [0,T^*]$. 
Since, by assumption \eqref{ass:thm:4},
\[
   4\tau < \upsi^0 = (N-1) \min_{i,j\in [N]}\psi_{ij} (0),
\]
and the weights $\psi_{ij}=\psi_{ij}(t)$ depend continuously on the group diameter
in the sense of \eqref{ass:psi:cont}, which in turn depends continuously on time
along the solutions of \eqref{eq:HKreact},
there exists $\eps>0$ and $\upsi > 4\tau$ such that
\[
    \upsi \leq (N-1) \min_{i,j\in [N]} \psi_{ij} (t)  \qquad \mbox{for } t\in [0, T^*+\eps].
\]
Lemma \ref{lem:4} gives then
\[
   \tot{}{t} d_x (t) \leq 4 \int_{t-\tau}^t d_x(s-\tau) \d s - \upsi d_x (t),
\]
and Corollary \ref{cor:Halanay:2} with $\alpha:=4\tau$ and $\beta:=\upsi$
yields
\[
   d_x(t) \leq d_x^0 e^{-Ct} \leq d_x^0 \qquad \mbox{for } t\in [0, T^*+\eps]
\]
with $C=C(\upsi)>0$ given by \eqref{Halanay:2:2}. This is a contradiction
to the choice of $T^*$ being the supremum of $\mathcal{S}$.

Therefore, we have $d_x(t) \leq d_x^0$ for all $t\geq 0$ and, consequently, by \eqref{ass:psi:cont},
\[
    \min_{i,j\in [N]} \psi_{ij} (t) \geq \frac{\upsi^0}{N-1} \qquad \mbox{for } t\geq 0.
\]
An application of Lemma \ref{lem:4} and Corollary \ref{cor:Halanay:2} with $\alpha:=4\tau$ and $\beta:=\upsi^0$
finally gives \eqref{concl:thm:React2}.
\end{proof}


\end{document}